\documentclass[12pt,letterpaper]{amsart}
\makeatletter
\@namedef{subjclassname@2020}{%
\textup{2020} Mathematics Subject Classification}
\makeatother
\usepackage{geometry,xcolor,graphicx}
\usepackage{amsmath,amssymb,amsfonts,amsthm,amscd}
\usepackage{enumerate}
\usepackage{mathrsfs,latexsym,url,setspace}
\usepackage{mathtools}
\usepackage{mathpazo} 
\usepackage[colorlinks,allcolors=blue]{hyperref} 
\usepackage{xpatch} 
\xpatchcmd{\proof} 
{\itshape} 
{\bfseries} 
{}
{}
\newtheorem{theorem}{Theorem}
\newtheorem{proposition}{Proposition}
\newtheorem{lemma}{Lemma}
\newtheorem{corollary}{Corollary}
\theoremstyle{remark}
\newtheorem{remark}{Remark}
\newtheorem{example}{Example}
\newcommand{\C}{\mathbb{C}}

\newcommand{\disk}{\mathbb{D}}
\newcommand{\T}{\mathbb{T}}
\newcommand{\ball}{\mathbb{B}}
\newcommand{\ballb}{\overline{\mathbb{B}}}

\newcommand{\zb}{\overline{z}}
\newcommand{\wb}{\overline{w}}
\allowdisplaybreaks
\onehalfspace

\title{On compactness of products of Toeplitz operators}

\author{Trieu Le}
\address{University of Toledo, Department of 
	Mathematics \& Statistics, 2801 W. Bancroft, Toledo, OH 43606, USA}
\email{trieu.le2@utoledo.edu}
	
\author{Tomas Miguel Rodriguez}
\address{Mathematics and Statistics, Kenyon College, 106 College Park Dr, Gambier, OH 43022, USA}
\email{rodriguez2@kenyon.edu}
	
\author{S\"{o}nmez \c{S}ahuto\u{g}lu}
\address{University of Toledo, Department of 
	Mathematics \& Statistics, 2801 W. Bancroft, Toledo, OH 43606, USA}
\email{sonmez.sahutoglu@utoledo.edu}

\date{\today}
\begin{document}

\begin{abstract}
We study compactness of product of Toeplitz operators with symbols continuous 
on the closure of the polydisc in terms of behavior of the symbols on the boundary. 
For certain classes of symbols $f$ and $g$, we show that $T_fT_g$ is compact if and only 
if $fg$ vanishes on the boundary. We provide examples to show that for more general 
symbols, the vanishing of $fg$ on the whole polydisc might not imply the compactness 
of $T_fT_g$. On the other hand, the reverse direction is closely related to the zero 
product problem for Toeplitz operators on the unit disc, which is still open.
\end{abstract}

\maketitle

Let $\Omega$ be a bounded domain in $\C^n$. The Bergman space $A^2(\Omega)$ 
consists of all holomorphic functions on $\Omega$ that are square integrable 
with respect to the Lebesgue volume measure $dV$. The orthogonal projection 
$P:L^2(\Omega)\to A^2(\Omega)$ is known as the Bergman projection. 
For a bounded measurable function $f$ on $\Omega$, the Toeplitz operator 
$T_{f}: A^2(\Omega)\rightarrow A^2(\Omega)$ is defined as
\[T_fh = P(fh)\] 
for $h\in A^2(\Omega)$. We call $f$ the symbol of $T_f$. 

There is an extensive literature on the study of Toeplitz operators on various 
domains. In this paper, we are particularly interested in the case the domain 
is the polydisc and compactness of product of Toeplitz operators whose symbols 
are continuous up to the boundary. 

A classical approach to compactness of Toeplitz operators involves the Berezin 
transform. For finite sum of finite products of Toeplitz operators on the 
Bergman space of the unit disc, the Axler--Zheng Theorem 
\cite[Theorem 2.2]{AxlerZheng98} characterizes compactness in terms of the 
behavior of the Berezin transform of the operator. In higher dimensions, 
the Axler--Zheng Theorem is extended to the case of the polydisc as seen 
in \cite{Englis1999} and \cite[p. 232]{ChoeKooLee09}, and the unit ball as shown in 
\cite[Theorem 9.5]{Suarez07}. Recently, there have been a few generalizations 
of this result in different directions. See, for instance, 
\cite{CuckovicSahutoglu2013,MitkovskiSuarezWick2013,MitkovskiWick2014,CuckovicSahutogluZeytuncu2018}.

In this paper, we study compactness of products of Toeplitz operators in terms of 
the behavior of the symbols on the boundary. More specifically, we would like to 
characterize functions $f,g$ that are continuous on $\overline{\disk^n}$ such that 
$T_{f}T_{g}$ is compact.

Coburn \cite[Lemma 2]{Coburn73} showed that on the Bergman space over unit 
ball $\ball$, for $f$ a continuous function on $\ballb$, the Toeplitz operator 
$T_f$ is compact if and only if $f=0$ on $b\ball$. 
Furthermore, \cite[Theorem 1]{Coburn73} established a $*$-isomorphism 
$\sigma:\tau(\ball) / \mathscr{K} \rightarrow C(b\ball)$ satisfying
\[\sigma(T_{f}+ \mathscr{K}) = f|_{b\ball},\]
where $\tau(\ball)$ is the Toeplitz algebra generated by 
$\{T_{\varphi}: \varphi \in C(\overline{\ball})\}$ and $\mathscr{K}$  is the ideal of 
compact operators on $A^2(\ball)$. As a consequence, we see that for 
$f_1,\ldots,f_N\in C(\ballb)$, the product $T_{f_1}\cdots T_{f_N}$ is compact 
if and only if the product $f_1\cdots f_N = 0$ on $b\ball$.

On the polydisc $\disk^n$, the first author \cite{Le2010} showed that, in the context 
of weighted Bergman spaces, for $f\in C(\overline{\disk^n})$, the Toeplitz 
operator $T_f$ is compact if and only if $f$ vanishes on $b\disk^n$, 
the (topological) boundary of $\disk^n$. Generalizing 
this result, the second and the third authors in \cite{RodriguezSahutogluPreprint} 
proved that compactness of the Toeplitz operator with a symbol continuous on the 
closure of a bounded pseudoconvex domain in $\C^n$ with Lipschitz boundary is 
equivalent to the symbol vanishing on the boundary of the domain. 

Motivated by Coburn's aforementioned result, one may expect that the necessary 
and sufficient condition for $T_{f}T_{g}$ to be compact is that $fg$ vanishes on 
$b\disk^n$. However, we shall see in our 
results and examples that while the above statement holds for a certain class of 
symbols, sufficiency is false in general. On the other hand, necessity is closely 
related with the famous ``zero product problem'' in the theory of Toeplitz 
operators on the unit disc, which is still wide open.

\section{Main Result}
Let $T= \sum^{N}_{j=1} T_{f_{j,1}} \cdots T_{f_{j,m_j}}$ be a finite sum of finite 
products of Toeplitz operators with $f_{j,k}\in C(\overline{\disk})$. 
Coburn's aforementioned result implies that compactness of $T$ on $A^2(\disk)$ 
is equivalent to $\sum^{N}_{j=1} {f_{j,1}} \cdots {f_{j,m_j}} = 0$ on the circle. 
Therefore, throughout the paper we will assume that $n\geq 2$ as  
the case $n=1$ is well understood.

Before we state our results, we define the restriction operator 
$R_{k,\xi}:C(\overline{\disk^n})\to  C(\overline{\disk^{n-1}})$ 
for $\xi\in \T$ and $k=1,\ldots,n$ as follows. 
\begin{align*}
R_{1,\xi}f(z_1,\ldots,z_{n-1})=&\,f(\xi,z_1,\ldots,z_{n-1}),\\
R_{n,\xi}f(z_1,\ldots,z_{n-1})=&\,f(z_1,\ldots,z_{n-1},\xi), 
\end{align*} 
and 
\[R_{k,\xi}f(z_1,\ldots,z_{n-1})=f(z_1,\ldots,z_{k-1},\xi,z_k,\ldots,z_{n-1})\]
for $2\leq k\leq n-1$ and $f\in C(\overline{\disk^n})$. 

In our main result, we give a characterization of compactness of the 
finite sum of finite products of Toeplitz operators in terms of the vanishing 
of the operator restricted to the polydiscs in the boundary. We recall that 
$b\disk^n$ consists of all $z=(z_1,\ldots,z_n)\in \overline{\disk^n}$ such 
that $|z_j|=1$ for some $j$. 

\begin{theorem}\label{ThmRestrictionSlice}
Let  $T= \sum^{N}_{j=1} T_{f_{j,1}} \cdots T_{f_{j,m_j}}$ be a finite sum of 
finite products of Toeplitz operators on $A^2(\disk^n)$ for 
$f_{j,k} \in C(\overline{\disk^n})$ with $n\geq 2$. Then $T$ is compact on 
$A^2(\disk^n)$ if and only if 
\[\sum^{N}_{j=1} T_{R_{k,\xi}f_{j,1}} \cdots T_{R_{k,\xi}f_{j,m_j}} = 0\] 
on $A^2(\disk^{n-1})$ for all $\xi\in\T$ and $1\leq k\leq n$.
\end{theorem}

As an immediate corollary we get the following. 

\begin{corollary}  \label{CorDirection1}
Let $f_j\in C(\overline{\disk^n})$ for $1\leq j\leq m$. Assume that for each 
$\xi\in\T$ and $1\leq k \leq n$ there exists $j$ such that  $R_{k,\xi}f_j=0$ 
on $\disk^{n-1}$. Then $T_{f_m}\cdots T_{f_1}$ is compact on $A^2(\disk^n)$.
\end{corollary}

\section{Applications}
Let $\varphi$ and $\psi$ be two functions in $C(\overline{\disk})$. 
We define $f(z,w) = \varphi(w)$ and $g(z,w) = \psi(w)$ for 
$z,w\in\overline{\disk}$. Then for any $\xi\in\T$,
\[R_{1,\xi}f(w) = f(\xi,w) = \varphi(w),\quad R_{1,\xi}g(w) 
= g(\xi,w)=\psi(w)\quad\text{ for } w\in\disk\] 
and
\[R_{2,\xi}f(z) = \varphi(\xi),\quad R_{2,\xi}g(z) 
= \psi(\xi)\quad\text{ for } z\in\disk.\]
By Theorem \ref{ThmRestrictionSlice}, the product $T_fT_g$ is 
compact on $A^2(\disk^2)$ if and only if $T_{\varphi}T_{\psi} = 0$ 
on $A^2(\disk)$ and $\varphi(\xi)\psi(\xi) = 0$ for all $\xi\in\T$. 
Since the second condition is actually a consequence of the first, 
we conclude that for such $f$ and $g$, the product $T_fT_g$ is 
compact on $A^2(\disk^2)$ if and only if $T_{\varphi}T_{\psi}= 0$ 
on $A^2(\disk)$, which is equivalent to $T_{f}T_{g} = 0$ on $A^2(\disk^2)$.

\begin{example}\label{E:fg_is_zero}
Let
\[\varphi(w) = \begin{cases}
1-2|w| & \text{ for } 0\leq |w|\leq\frac{1}{2}\\
0 & \text{ for } |w|>\frac{1}{2},
\end{cases}\]
and
\[\psi(w) = \begin{cases}
0 & \text{ for } 0\leq |w|\leq\frac{1}{2}\\
2|w|-1 & \text{ for } |w|>\frac{1}{2}.
\end{cases}\]
Using polar coordinates, one can check that both operators $T_{\varphi}$ 
and $T_{\psi}$ are diagonalizable 
with respect to the standard orthonormal basis and their eigenvalues are 
all strictly positive. Hence  $T_{\varphi}T_{\psi}\not\equiv 0$ on $A^2(\disk)$.
On the other hand, $\varphi\psi = 0$ on $\overline{\disk}$. Then 
for $f(z,w)=\varphi(w)$ and $g(z,w)=\psi(w)$, we have $fg = 0$ 
on $\overline{\disk^2}$ but $T_{f}T_{g}$ is not compact on $A^2(\disk^2)$ 
as $T_{\varphi}T_{\psi}\not\equiv 0$. This example  shows that the vanishing of $fg$ 
on $b\disk^2$ (or even on $\overline{\disk^2}$) does not imply the 
compactness of $T_{f}T_{g}$. 
\end{example}

\begin{example}\label{E:fg_not_zero}
Take $f$ as in Example \ref{E:fg_is_zero} and define
\[g(z,w) = \varphi(z) + \psi(w).\] 
Then $fg$ is not identically zero on $\disk^2$ because $f(0,0) = g(0,0) = 1$ 
and $fg = 0$ on $b\disk^2$. Yet, by Theorem \ref{ThmRestrictionSlice}, 
the product $T_fT_g$ is not compact since for $\xi\in\T$,
\[T_{R_{1,\xi}f}T_{R_{1,\xi}g} = T_{\varphi}T_{\psi}\] 
is not the zero operator on $A^2(\disk)$.
\end{example}

\begin{remark} 
\label{R:zero_product}
From the previous examples we see that $fg=0$ on $b\disk^2$ is not 
a sufficient condition for the compactness of $T_fT_g$. Is it a 
necessary condition? It turns out this question is  related to 
the \textit{zero product problem} for Toeplitz operators on the 
disc. More specifically, as in 
Example \ref{E:fg_is_zero}, we see that with $f(z,w)=\varphi(w)$ 
and $g(z,w) = \psi(w)$, if the product $T_fT_g$ is compact on 
$A^2(\disk^2)$, then $T_{\varphi}T_{\psi} = 0$ on $A^2(\disk)$ 
(which gives $\varphi\psi = 0$ on $\T$). However, it is not known 
if this condition implies that $\varphi\psi = 0$ on $\disk$.
For $\xi\in\T$ and $z, w\in\disk$, we have
$f(\xi,w)g(\xi,w) = \varphi(w)\psi(w)$ and $f(z,\xi)g(z,\xi) 
= \varphi(\xi)\psi(\xi)$. 
So $fg = 0$ on $b\disk^2$ if and only if $\varphi\psi = 0$ on $\disk$.
\end{remark}

In Proposition \ref{CorHarmonic} below, we show that if the symbols 
are harmonic along the discs in the boundary, then we have necessary 
and sufficient conditions for the compactness of the product of two 
Toeplitz operators. A function $f \in C^2(\disk^n)$ is said to be 
$n$-harmonic if 
\[\Delta_j f = 4\dfrac{\partial^2 f}{\partial z_j \partial \overline{z}_j}=0,\]
for all $j=1,2,\dots,n.$ That is, $f$ is harmonic in each variable 
separately \cite[pg. 16]{Rudin69}. 

\begin{proposition}\label{CorHarmonic}
Let $f,g\in C(\overline{\disk^n})$ (with $n\geq 2$) such that for 
$\xi\in\T$, and $1\leq k\leq n$, the functions $R_{k,\xi}f$ and  
$R_{k,\xi}g$ are $(n-1)$-harmonic on $\disk^{n-1}$. Then $T_fT_g$ is 
compact if and only if $fg = 0$ on $b\disk^n$.
\end{proposition}

We note that in Example \ref{E:fg_is_zero}, both $f$ and $g$ depend on the 
same single variable. In Proposition \ref{PropDecoupled} below, we give a 
characterization when the symbols are product of single-variable functions.

\begin{proposition} \label{PropDecoupled}
Let  $T=\prod^M_{k=1}T_{f_k}$ be a finite product of Toeplitz operators 
on $A^2(\disk^n)$ such that $f_k(z)=\prod^{n}_{j=1}f_{j,k}(z_j)$ for  
$f_{j,k} \in C(\overline{\disk})$ and $z=(z_1,\ldots,z_n)\in\disk^n$. 
Let $F=\prod^M_{k=1} f_k$. Then the following statements hold.
\begin{enumerate}[(i)]
\item If $T$ is a nonzero compact operator, then $F = 0$ on $b\disk^n$.
\item If $F = 0$ on $b\disk^n$ and $F$ is not identically zero on $\disk^n$, 
then $T$ is compact.
\end{enumerate}
\end{proposition}

\begin{remark}
We do not know whether (i) in Proposition \ref{PropDecoupled} still holds 
in the case $T$ is the zero operator. This is closely related to the zero 
product problem. More specifically, consider $f(z,w)=\varphi(w)$ and 
$g(z,w) = \psi(w)$, where $\varphi, \psi\in C(\overline{\disk})$. 
Then $T = T_fT_g = 0$ on $A^2(\disk^2)$ if and only if $T_{\varphi}T_{\psi} = 0$ 
on $A^2(\disk)$. On the other hand, $F = fg = 0$ on $b\disk^n$ 
if and only if $\varphi\psi = 0$ on $\disk$. It is still an open 
problem whether $T_{\varphi}T_{\psi} = 0$ on $A^2(\disk)$ 
implies that $\varphi\psi = 0$ on $\disk$.
\end{remark}

\begin{remark}
The conclusion of (ii) in Proposition \ref{PropDecoupled} does not hold if 
$F$ is identically zero on $\disk^n$. Indeed, the functions $f$ and $g$ 
in Example \ref{E:fg_is_zero} are of the type considered here and 
$F = fg = 0$ on $\overline{\disk^2}$ but $T_{f}T_{g}$ is not compact 
on $A^2(\disk^2)$.
\end{remark}

In the proposition below, we show that when all but at most one of the 
symbols are polynomials, compactness of a Toeplitz product on 
$A^2(\disk^2)$ is equivalent to the vanishing of the product 
of the symbols on $b\disk^2$. For this result, we need to 
restrict to dimension two. It would be interesting to extend 
the result to all $n\geq 2$. See Remark \ref{R:ThmPolynomial}.

\begin{proposition}\label{ThmPolynomial}
Let $f_1,\ldots,f_M$ and $g_1,\ldots,g_N$ be polynomials in $z,w$ and 
$\zb,\wb$, and $h\in C(\overline{\disk^2})$. 
Then $T_{f_1}\cdots T_{f_M}T_hT_{g_1}\cdots T_{g_N}$ is compact on 
$A^2(\disk^2)$ if and only if 
\[f_1\cdots f_M h g_1\cdots g_N = 0 \text{ on } b\disk^2.\]
\end{proposition}

\section{Proofs}
Let $BT(p)$ denote the Berezin transform of a bounded linear operator 
$T:A^2(\disk^n)\to A^2(\disk^n)$ at $p\in \disk^n$. That is, 
\[BT(p)=\langle Tk_p,k_p\rangle\]
where 
\[k_p(z)=\frac{K(z,p)}{\sqrt{K(p,p)}}\] 
is the normalized Bergman kernel of $\disk^n$.

We will need the following lemma whose proof is contained in the proof of 
Theorem 1 in \cite{CuckovicHuoSahutogluPreprint}. We provide a sketch of
the proof here for the convenience of the reader. We note that $Bf$ denotes 
$BT_f$ whenever $f$ is a bounded function and we use the following notation: 
$z' =(z_2,\ldots,z_n)\in \C^{n-1}$ for $z=(z_1,\ldots,z_n)\in \C^n$. 
For functions $h_1$ defined on $\disk$ and $h_2$ defined on $\disk^{n-1}$, 
we use $h_1h_2$ to denote the function $h_1(z_1)h_2(z')$ on $\disk^n$.

\begin{lemma}\label{Lem1} 
Suppose $n\geq 2$ and $\psi\in C(\overline{\disk^n})$. 
Let $q=(\zeta,q')\in \T\times\overline{\disk^{n-1}}$ and define 
$\psi_{\zeta}(z) = \psi(\zeta,z')$ for $z\in\disk^n$. 
\begin{enumerate}[(i)]
\item If $\{h_p: p\in\disk^{n}\}$ is a bounded set in $L^2(\disk^{n-1})$, then
\[\lim_{p\to q}\big\|(\psi - \psi_{\zeta})k_{p_1}^{\disk}h_p\big\| = 0.\]
\item If $\psi_1,\ldots,\psi_v\in C(\overline{\disk^{n}})$ are functions 
independent of $z_1$ and $W$ is any bounded operator on $L^2(\disk^{n})$, then
\[\lim_{p\to q}\big\|WT_{\psi-\psi_{\zeta}}T_{\psi_1}\cdots T_{\psi_v}k_p\big\| = 0.\]
\end{enumerate}
\end{lemma} 

\begin{proof}
(i) 
Let $\epsilon>0$ be given. By the uniform continuity of $\psi$, 
there exists $\delta > 0$ such that for all $z'\in\disk^{n-1}$,
\[|\psi(z_1,z')-\psi_{\xi}(z_1,z')| 
< \dfrac{\epsilon}{\sup\{\|h_p\|_{L^2(\disk^{n-1})}\}+1} 
\text{ whenever }\ |z_1 - \xi| < \delta. \]
Then,
\begin{align*}
\|(\psi - \psi_{\xi})k^{\disk}_{p_1}h_p\|^2 
& = \|(\psi - \psi_{\xi})k^{\disk}_{p_1}h_p\|^2_{L^2( \{z \in \disk^n:|z_1 - \xi| 
< \delta\})}\\
& \qquad\qquad+ \|(\psi - \psi_{\xi})k^{\disk}_{p_1}h_p\|^2_{L^2( \{z \in\disk^n:
|z_1 - \xi| \geq \delta\})} \\
& \leq \epsilon^2 + \pi\|h_p\|^2_{L^2(\disk^{n-1})}\|(\psi 
- \psi_{\xi})k^{\disk}_{p_1}\|^2_{L^{\infty}(\{z \in \disk^n:|z_1 - \xi| 
\geq \delta\})}. 
\end{align*}
However, 
\[\sup\left\{\left|k^{\disk}_{p_1}(z_1)\right|
: |z_1 - \xi| \geq \delta\right\} \rightarrow 0\ \text{as}\ p_1 \rightarrow \xi.\] 
Then, $\limsup_{p \rightarrow q}\|(\psi -\psi_{\xi})k^{\disk}_{p_1}h_p\| \leq \epsilon$. 
Since $\epsilon>0$ was arbitrary, we conclude that
\[\lim_{p\to q}\big\|(\psi - \psi_{\zeta})k_{p_1}^{\disk}h_p\big\| = 0.\]

(ii) We note that $k_p = k_{p_1}^{\disk}k_{p'}^{\disk^{n-1}}$ for $p = (p_1,p')$. 
We define
\[h_p = T_{\psi_1}\cdots T_{\psi_v}k_{p'}^{\disk^{n-1}}\text{ for } p\in\disk^{n}.\] 
Since each $\psi_j$ is independent of $z_1$, $h_p$ is independent of $z_1$ and 
hence it can be considered as an element of $L^2(\disk^{n-1})$. Note that the 
set $\{h_p: p\in\disk^{n}\}$ is bounded by $\|T_{\psi_1}\cdots T_{\psi_v}\|$. 
Furthermore, we have
$T_{\psi_1}\cdots T_{\psi_v}k_p = k_{p_1}^{\disk}h_p$. It follows that
\[\big\|WT_{\psi-\psi_{\zeta}}T_{\psi_1}\cdots T_{\psi_v}k_p\big\|
\leq \|W\|\cdot\big\|(\psi-\psi_{\zeta})k_{p_1}^{\disk}h_p\big\|,\] which, by (i), 
converges to zero as $p\to q$.
\end{proof}

\begin{proof}[Proof of Theorem \ref{ThmRestrictionSlice}]
We first make an observation. If $\varphi$ is a bounded function on 
$\disk^{n-1}$, then $T_{\varphi}$, while initially defined on $A^2(\disk^{n-1})$, 
can be naturally considered as a Toeplitz operator with 
symbol $E_1\varphi(z_1,z') =\varphi(z')$ acting on $A^2(\disk^n)$. 
This will not create any confusion 
due to the fact that for $h\in A^2(\disk^n)$ independent of $z_1$, the 
function $T_{E_1\varphi}h$ is also independent of $z_1$ and 
$(T_{E_1\varphi}h)(z)=(T_{\varphi}h)(z')$ for all $z=(z_1,z')\in\disk^{n}$.

Let $\xi\in\T$. For each $j$ and $m_j$, the function $f_{j,m_j}$ 
can be written as $f_{j,m_j} = (f_{j,m_j} -R_{1,\xi}f_{j,m_j})+R_{1,\xi}f_{j,m_j}$.
We then expand $T=\sum^{N}_{j=1} T_{f_{j,1}} \cdots T_{f_{j,m_j}}$ as 
\begin{align*} 
& T = \sum^{N}_{j=1} \left(T_{R_{1,\xi}f_{j,1}} \cdots T_{R_{1,\xi}f_{j,m_j}} 
+ T_{f_{j,1} -R_{1,\xi}f_{j,1}}T_{R_{1,\xi}f_{j,2}}
\cdots T_{R_{1,\xi}f_{j,m_j}} \right. \\ 
& \left. + T_{f_{j,1}}T_{f_{j,2} - R_{1,\xi}f_{j,2}}T_{R_{1,\xi}f_{j,3}} 
\cdots T_{R_{1,\xi}f_{j,m_j}} +\cdots+T_{f_{j,1}}T_{f_{j,2}} 
\cdots T_{f_{j,m_j-1}}T_{f_{j,m_j}-R_{1,\xi}f_{j,m_j}}\right)  \\
& = \sum^{N}_{j=1} T_{R_{1,\xi}f_{j,1}} \cdots T_{R_{1,\xi}f_{j,m_j}} 
+ \sum^{N}_{j=1}\left(T_{f_{j,1} -R_{1,\xi}f_{j,1}}T_{R_{1,\xi}f_{j,2}}
\cdots T_{R_{1,\xi}f_{j,m_j}} \right. \\ 
& \left. + T_{f_{j,1}}T_{f_{j,2} - R_{1,\xi}f_{j,2}}T_{R_{1,\xi}f_{j,3}} 
\cdots T_{R_{1,\xi}f_{j,m_j}} +\cdots+T_{f_{j,1}}T_{f_{j,2}} 
\cdots T_{f_{j,m_j-1}}T_{f_{j,m_j}-R_{1,\xi}f_{j,m_j}}\right).
\end{align*}
Note that in the second sum, each summand has the form considered in 
Lemma \ref{Lem1}(ii). We then conclude that for any 
$q=(\xi,q')\in\T\times\overline{\disk^{n-1}}$,
\begin{align}\label{Eqn:limitTk_p}
\lim_{p\to q}\Big\|Tk_p - \sum^{N}_{j=1} T_{R_{1,\xi}f_{j,1}} 
\cdots T_{R_{1,\xi}f_{j,m_j}}k_p\Big\| = 0.
\end{align}
Now suppose that $T$ is compact. Fix $p'\in\disk^{n-1}$. Since $k_{(p_1,p')}\to 0$ 
weakly as $p_1\to\xi$, the compactness of $T$ implies that $\|Tk_{(p_1,p')}\|\to 0$ 
as $p_1 \to \xi$. Equation \eqref{Eqn:limitTk_p} then gives
\begin{align}\label{EqnLim}
\lim_{p_1\to\xi}\Big\|\sum^{N}_{j=1} T_{R_{1,\xi}f_{j,1}} 
\cdots T_{R_{1,\xi}f_{j,m_j}}k_{(p_1,p')}\Big\| = 0.    
\end{align}
Since
\begin{align*}
\sum^{N}_{j=1} T_{R_{1,\xi}f_{j,1}} \cdots T_{R_{1,\xi}f_{j,m_j}}k_{(p_1,p')} 
& = \sum^{N}_{j=1} T_{R_{1,\xi}f_{j,1}} 
\cdots T_{R_{1,\xi}f_{j,m_j}}(k_{p_1}^{\disk}k_{p'}^{\disk^{n-1}}) \\
& = k^{\disk}_{p_1}\cdot\sum^{N}_{j=1} T_{R_{1,\xi}f_{j,1}} 
\cdots T_{R_{1,\xi}f_{j,m_j}}k^{\disk^{n-1}}_{p'}
\end{align*}
and $\|k_{p_1}^{\disk}\|=1$ for all $p_1$, \eqref{EqnLim} implies that
\[\sum^{N}_{j=1} T_{R_{1,\xi}f_{j,1}} 
\cdots T_{R_{1,\xi}f_{j,m_j}}k^{\disk^{n-1}}_{p'} = 0.\]
Because $p'$ was arbitrary, it follows that 
$\sum^{N}_{j=1}  T_{R_{1,\xi}f_{j,1}} \cdots T_{R_{1,\xi}f_{j,m_j}}$ 
is the zero operator on $A^2(\disk^{n-1})$. Applying the same method 
for other values of $k$, we have 
\[\sum^{N}_{j=1} T_{R_{k,\xi}f_{j,1}} \cdots T_{R_{k,\xi}f_{j,m_j}} = 0\]
on $A^2(\disk^{n-1})$ for  $1 \leq k\leq n$ and all $\xi \in \T$.

Let us now prove the converse. Let $q=(\xi,q') \in b\disk^n$ 
with $\xi \in \T$ and $q'\in\overline{\disk^{n-1}}$. Since 
$\sum^{N}_{j=1}  T_{R_{1,\xi}f_{j,1}} \cdots T_{R_{1,\xi}f_{j,m_j}} = 0$, 
equation \eqref{Eqn:limitTk_p} implies that $\lim_{p\to q}\|Tk_p\| = 0$. 
As a consequence, 
\[\lim_{p\to q}BT(p) = \lim_{p\to q}\langle Tk_p,k_p\rangle = 0.\] 
The same argument is applicable for all $q\in b\disk^{n}$. By Axler--Zheng Theorem 
for $\disk^n$ (\cite{Englis1999} and \cite[p. 232]{ChoeKooLee09}), we conclude that 
$T$ is compact on $A^2(\disk^n)$.
\end{proof}

\begin{proof}[Proof of Corollary \ref{CorDirection1}]
We assume that for each $\xi\in\T$ and $1\leq k \leq n$ there exists $j$ 
such that $R_{k,\xi}f_j=0$. Then $T_{R_{k,\xi}f_m}\cdots T_{R_{k,\xi}f_{1}}=0$ 
on $A^2(\disk^{n-1})$. Hence, Theorem \ref{ThmRestrictionSlice} implies that 
$T_{f_m}\cdots T_{f_1}$ is compact on $A^2(\disk^n)$.		
\end{proof} 

\begin{proof}[Proof of Proposition \ref{CorHarmonic}] 
To prove the forward direction, we first use Theorem \ref{ThmRestrictionSlice} 
to conclude that $T_{R_{k,\xi}g}T_{R_{k,\xi}f}$ 
is the zero operator on $A^2(\disk^{n-1})$ for all $\xi\in \T$ and 
$1\leq k \leq n$. Since the symbols $R_{k,\xi}f$ and $R_{k,\xi}g$ 
are $(n-1)$-harmonic on $\disk^{n-1}$, we apply \cite[Theorem 1.1]{ChoeKooLee07} 
(or \cite[Corollary 2]{AhernCuckovic2001} in the case $n=2$) to conclude that 
either $R_{k,\xi}f=0$ or $R_{k,\xi}g=0$. Then $fg=0$ on $b\disk^n$ as desired.

To prove the converse we argue as follows. For each $1\leq k\leq n$ and $\xi\in \T$, 
since both $R_{k,\xi}f$ and $R_{k,\xi}g$ are $(n-1)$-harmonic and their product is 
zero on $\disk^{n-1}$, either $R_{k,\xi}f=0$ or $R_{k,\xi}g=0$. Then 
$T_{R_{k,\xi}g}T_{R_{k,\xi}f}=0$ on $A^2(\disk^{n-1})$ for all $\xi \in \T$ and 
$1 \leq k \leq n$. Theorem \ref{ThmRestrictionSlice} now implies that $T_gT_f$ 
is compact.
\end{proof} 

\begin{proof}[Proof of Proposition \ref{PropDecoupled}] 
We first prove (i). Assume that $T$ is a nonzero compact operator. 
Then by Theorem \ref{ThmRestrictionSlice} when restricted on the 
first coordinate, for any $\xi\in\T$,
\[0 = \prod_{k=1}^{M}T_{R_{1,\xi}f_k}
=\Big(\prod_{k=1}^{M}f_{1,k}(\xi)\Big)\prod_{k=1}^{M}T_{\widetilde{f}_k}\] 
on $A^2(\disk^{n-1})$, where $\widetilde{f}_k(z_2,\ldots,z_n) = f_{2,k}(z_2)\cdots f_{n,k}(z_n)$. 
Since $T$ is not the zero operator, the second factor on the right hand side above 
is a nonzero operator. This follows from the fact that $T$ can be written as 
the product 
\[\displaystyle\big(\prod_{k=1}^{M}T_{f_{1,k}}\big)
\cdot\big(\prod_{k=1}^{M}T_{\widetilde{f}_k}\big)\] 
where the first factor acts on functions in $z_1$ and the second factor 
acts on functions in $z' = (z_2,\ldots,z_n)$. Hence, 
$\prod_{k=1}^{M}f_{1,k}(\xi)=0$. It follows that
\[F(\xi,z_2,\dots,z_n)=\prod_{k=1}^{M}f_k(\xi,z_2,\dots,z_n)
=\left(\prod_{k=1}^{M}f_{1,k}(\xi)\right)\left(\prod^{n}_{j=2}
\prod_{k=1}^{M}f_{j,k}(z_j)\right)=0\]
on $\T \times \disk^{n-1}$. The same argument applies to other 
coordinates and we have $F = 0$ on $b\disk^n$.

Next we prove (ii). Assume that $F=\prod^M_{k=1} f_k = 0$ on $b\disk^n$ and $F$ is 
not identically zero on $\disk^n$. Choose $q=(q_1,\ldots,q_n)\in\disk^n$ 
such that $f_k(q)\neq 0$ for all $k$, which implies that $f_{j,k}(q_j)\neq 0$ 
for all $j$ and $k$. For any $\xi\in\T$, since $z=(\xi,q_2,\ldots,q_n)\in b\disk^n$, 
we have
\[0 = F(z) = \Big(\prod_{k=1}^{M}f_{1,k}(\xi)\Big)
\cdot\prod_{j=2}^{n}\prod_{k=1}^{M}f_{j,k}(q_j).\]
Because the second factor is nonzero, it follows that 
$\prod_{k=1}^{M}f_{1,k}(\xi) = 0$. As a result,
\[\prod_{k=1}^{M}T_{R_{1,\xi}f_k} 
= \Big(\prod_{k=1}^{M}f_{1,k}(\xi)\Big)\prod_{k=1}^{M}T_{\widetilde{f}_k} = 0\] 
on $A^2(\disk^{n-1})$, where, as before, 
$\widetilde{f}_k(z_2,\ldots,z_n) = f_{2,k}(z_2)\cdots f_{n,k}(z_n)$. 
The same argument applies to other parts of $b\disk^n$. 
Then Theorem \ref{ThmRestrictionSlice} implies that $T=\prod_{k=1}^{M}T_{f_k}$ 
is compact on $A^2(\disk^n)$.
\end{proof}

The proof of Proposition \ref{ThmPolynomial} hinges on several elementary facts 
about polynomials that we describe below. We use $\C[z,\zb]$ to denote the 
vector space of all polynomials in $z$ and $\zb$. 

The following lemma is well known.  The proof follows from the fact that 
if a real analytic function vanishes on a non-empty open set, 
it must be identically zero. 

\begin{lemma}\label{L:zero_poly}
Let $f\in\C[z,\zb]$ be not identically zero. Then the set 
\[\{z\in\C: f(z) = 0\}\] 
has an empty interior.
\end{lemma}

\begin{lemma}\label{L:zero_poly_circle}
Let $f\in\C[z,\zb]$. Assume that there exist infinitely many $\xi\in\T$ 
such that $f(\xi)=0$. Then there is a polynomial $g\in\C[z,\zb]$ such that 
$f(z) = (1-|z|^2)g(z)$. In particular, $f(\xi)=0$ for all $\xi\in\T$.
\end{lemma}

\begin{proof}
For non-negative integers $s, t$, we write
\[\zb^s z^t = \begin{cases}
    |z|^{2s}\,z^{t-s} & \text{ if } t\geq s,\\
    |z|^{2t}\,\zb^{s-t} & \text{ if } t < s.
\end{cases}\]
As a result, there are integers $m,M\geq 0$ and polynomials $p_j$ (for $0\leq j\leq M$) 
and $q_j$ (for $0\leq j\leq m$) of a single variable such that
\[f(z) = \sum_{j=0}^{M}p_j(|z|^2)z^{j} + \sum_{j=0}^{m}q_j(|z|^2)\zb^{j}.\]
By the hypothesis, there exists infinitely many 
$\xi\in\T$ such that 
\[\sum_{j=0}^{M}p_j(1)\xi^{j} + \sum_{j=0}^{m}q_j(1)\overline{\xi}^j = f(\xi) =  0.\]
This implies that $p_j(1) = q_j(1) = 0$ for each $j$. As a consequence, all $p_j(r)$ 
and $q_j(r)$ are divisible by $1-r$. We then conclude that $f(z)$ is divisible by 
$1-|z|^2$, from which the conclusion of the lemma follows.
\end{proof}

\begin{lemma}\label{LemSymbolVanish}
Let $f(z,w)$ be a polynomial in $z,w,\zb,\wb$ and let $h \in C(\overline{\disk^2})$. 
Assume that $fh = 0$ on $b\disk^2$. Then $f|_{\T\times\overline{\disk}}=0$ 
or $h|_{\T\times\overline{\disk}}=0$ and $f|_{\overline{\disk}\times\T}=0$ 
or $h|_{\overline{\disk}\times\T}=0$.
\end{lemma}
\begin{proof}
Assume that $h$ does not vanish identically on $\T\times\overline{\disk}$. 
By continuity, there exist a non-empty arc $J\subseteq\T$ and a non-empty 
open set $V\subseteq\disk$ such that $h(\xi,w)\neq 0$ for all $\xi\in J$ 
and $w\in V$. It follows that $f(\xi,w)=0$ for all such $\xi$ and $w$. 
For each $\xi\in J$, applying Lemma \ref{L:zero_poly}, we conclude that 
$f(\xi,w)=0$ for all $w\in\overline{\disk}$. Then for each $w\in\overline{\disk}$, 
since $f(\xi,w)$ vanishes on $J$ (which is an infinite set), 
Lemma \ref{L:zero_poly_circle} implies that $f(\xi,w)=0$ 
for all $\xi\in\T$. Therefore, $f$ vanishes identically on 
$\T\times\overline{\disk}$. The proof for $\overline{\disk}\times\T$ is similar.
\end{proof}

\begin{lemma}[{\cite[Corollary 1.8]{ThilakarathnaThesis}}]
\label{L:zeroProduct}
Suppose $\varphi_1,\ldots,\varphi_M$ and $\psi_1,\ldots,\psi_N$ are polynomials 
of $z, \zb$ in $\disk$ and $g\in L^{2}(\disk)$. If 
$T_{\varphi_1}\cdots T_{\varphi_M}T_{g}T_{\psi_1}\cdots T_{\psi_N}=0$ 
on $A^2(\disk)$, then one of the symbols must be zero.
\end{lemma}

\begin{proof}[Proof of Proposition \ref{ThmPolynomial}]
Assume that $T_{f_1}\cdots T_{f_M}T_hT_{g_1}\cdots T_{g_N}$ is compact on $A^2(\disk^2)$, 
then by Theorem \ref{ThmRestrictionSlice}, 
\[ T_{R_{1,\xi}f_1}\cdots T_{R_{1,\xi}f_M}T_{R_{1,\xi}h}T_{R_{1,\xi}g_1}
\cdots T_{R_{1,\xi}g_N}
= 0\] 
on $A^2(\disk)$ for all $\xi\in\T$. By Lemma \ref{L:zeroProduct}, one of 
$R_{1,\xi}f_1,\ldots, R_{1,\xi}f_M$, $R_{1,\xi}h$, and $R_{1,\xi}g_1,\ldots,R_{1,\xi}g_N$ 
is a zero function on $\disk$. Thus, $f_1\cdots f_M h g_1\cdots g_N = 0$ on 
$\T \times \overline{\disk}$. Similar argument works for $\overline{\disk} \times \T$. 
Therefore, $f_1\cdots f_M h g_1\cdots g_N = 0$ on $b\disk^2$.

For the converse, by Lemma \ref{LemSymbolVanish}, one of the symbols is identically 
zero on $\T\times\overline{\disk}$. It then follows that 
\[T_{R_{1,\xi}f_1}\cdots T_{R_{1,\xi}f_M}T_{R_{1,\xi}h}T_{R_{1,\xi}g_1}
\cdots T_{R_{1,\xi}g_N}=0.\]
Similarly,  
\[T_{R_{2,\xi}f_1}\cdots T_{R_{2,\xi}f_M}T_{R_{2,\xi}h}T_{R_{2,\xi}g_1}
\cdots T_{R_{2,\xi}g_N}=0.\]
Therefore, by Theorem \ref{ThmRestrictionSlice}, we conclude that 
$T_{f_1}\cdots T_{f_M}T_hT_{g_1}\cdots T_{g_N}$ is compact on $A^2(\disk^2)$.  
\end{proof}

\begin{remark}\label{R:ThmPolynomial}
It is desirable to generalize Proposition \ref{ThmPolynomial} to $\disk^n$ for all 
$n\geq 2$. While Lemmas \ref{L:zero_poly}, \ref{L:zero_poly_circle} and 
\ref{LemSymbolVanish} remain true for all $n$, Lemma \ref{L:zeroProduct} 
has only been known for the disc. In order to extend Proposition \ref{ThmPolynomial} 
to all $n\geq 2$, one needs to prove a several-variable version of 
Lemma \ref{L:zeroProduct}. Some partial results have been obtained in the 
literature. For example, the main results of \cite{CuckovicHuoSahutoglu2022} 
imply that Lemma \ref{L:zeroProduct} holds in several variables when $g=1$ 
or when all $\varphi_j, \psi_k$ are monomials. As a result, 
Proposition \ref{ThmPolynomial} holds on $\disk^n$ for all $n\geq 2$ in 
the case $h=1$, or in the case all $f_j$ and $g_k$ are monomials.
\end{remark}


\section*{Acknowledgment}
Trieu Le is partially supported by Simons Foundation Travel Support for 
Mathematicians MPS-TSM-00002303.

\end{document}